\documentclass[11pt]{article}
\usepackage{amsmath,latexsym,amsthm}
\usepackage[english]{babel}
\usepackage{graphicx}
\usepackage{graphics}
\usepackage{amssymb}
\usepackage{color}
\usepackage{textcomp}
\usepackage{amsfonts}
\usepackage{longtable}
\usepackage{fancyhdr}
\usepackage{url}
\makeatletter

\newtheorem{theorem}{Theorem}
\newtheorem{lemma}{Lemma}

\begin{document}
\sloppy

\begin{center}
\textbf{Cauchy problem for the time-fractional Airy-type equation on the star-graph with edge-dependent order of derivative}\\
\end{center}

\begin{center}

{\sc  Sobirov~Zarif$^1$, Rakhimov~Kamoladdin$^2$}\\
{\it National University of Uzbekistan $^1$,\\
V.I. Romanovskiy Institute of Mathematics $^1$,\\
University of Exact and Social Sciences in Tashkent $^2$\\
}
e-mail: {\tt sobirovzar@gmail.com $^1$, \\
kamoliddin.ru1@gmail.com $^2$}

\end{center}
{\bf Abstract.}   We study the Cauchy problem for a time-fractional Airy-type equation on an infinite metric star graph with edge-dependent Caputo  derivatives. By constructing layer potentials generated by fractional Airy fundamental solutions, the problem is reduced to a generalized Abel integral equation. Using Pskhu’s solvability results, we establish existence, uniqueness, and an explicit integral representation of the solution. The obtained framework extends previous results where identical fractional orders were assumed on all edges.



\section{Introduction}

 We consider the star-graph $\Gamma$ with $k$ incoming bonds $b_i, i=\overline{1,k},$ and $m$ outgoing bonds $b_j, j=\overline{k+1,k+m}$.  The coordinates on the incoming edges are defined from $-\infty$ to $0$, and on the outgoing bonds from $0$ to $+\infty$. Let $u:\Gamma\to \mathbb{R}$ denotes a function on the whole graph. For each bond $b_n$ we define restriction of this function as $u_n=u|_{b_n}$. For the functions defined on the graph, we also use vector-type notations $u=(u_1,...,u_n)^T$, $u_x=\left(u_{1,x},..., u_{n,x}\right)^T, n=\overline{1,k+m}$. 

 On the star-graph  $\Gamma$ we consider the following problem

\begin{equation}\label{masala}
\begin{cases}
  \partial_{0t}^{\alpha_n} u_n(x,t) - u_{n,xxx}(x,t) = f_n(x,t), \,\,\,\ x\in b_n, \,\,(n=\overline{1,k+m}\,\ t\in (0,T],\\
  u_n(x,0) = u_{0,n}(x), \,\,\,\ x\in b_n, \,\,\,(n=\overline{1,k+m}\\
  A u(0,t) = 0, \,\,\,\ t\in[0,T],\\
  u_x^{+}(0,t) = B u_x^{-}(0,t), \,\,\,\ t\in[0,T],\\
  C^{-} \dfrac{\partial^2}{\partial x^2} u^{-}(x,t)\Big|_{x=0} = C^{+} \dfrac{\partial^2}{\partial x^2} u^{+}(x,t)\Big|_{x=0}, \,\,\,\ t\in[0,T].
\end{cases}
\end{equation}
Here $0<\alpha_n<1, \,\,\,\ (n=\overline{1,k+m}),$ $\,{{u}^{-}}={{\left( {{u}_{1}},{{u}_{2}},...,{{u}_{k}} \right)}^{T}},\,\,{{u}^{+}}={{\left( {{u}_{k+1}},{{u}_{k+2}},...,{{u}_{k+m}} \right)}^{T}},$ $A$ is the constant matrix of the form \[A=\left( \begin{matrix}
   0 & 0 & 0 & ... & 0  \\
   1 & -{{a}_{2}} & 0 & ... & 0  \\
   1 & 0 & -{{a}_{3}} & ... & 0  \\
   ... & ... & ... & ... & ...  \\
   1 & 0 & 0 & ... & 0  \\
   1 & 0 & 0 & ... & -{{a}_{k+m}}  \\
\end{matrix} \right),\]
${{C}^{-}}=(\frac{1}{{{a}_{1}}},\frac{1}{{{a}_{2}}},...,\frac{1}{{{a}_{k}}}),\,\,\,{{C}^{+}}=(\frac{1}{{{a}_{k+1}}},...,\frac{1}{{{a}_{k+m}}}),\,\,\,{{a}_{1}}=1,\,\,\,{{a}_{n}}\ne 0 \,\,\,\ (n=\overline{1,k+m})$ and $B$ is an $m\times k$ matrix of the form \[B=\left( \begin{matrix}
   {{b}_{k+1,1}} & ... & {{b}_{k+1,k}}  \\
   ... & ... & ...  \\
   {{b}_{k+m,1}} & ... & {{b}_{k+m,k}}  \\
\end{matrix} \right).\]
Here $\partial_{0t}^{\rho}g(t)$ stands for the Caputo fractional derivative
 \begin{equation} \label{KThos}
\partial_{0t}^{\rho }g(t)=\frac{1}{\Gamma(1-\rho)}\int\limits_{0}^{t}\frac{g'(\xi)}{(t-\xi)^\rho}d\xi, \ \,\,\,\ 0<\rho<1,
\end{equation}
where $\Gamma(\rho) \,\,\, - $ is the gamma function.

We suppose that the functions $u_n(x,t)$ and their derivatives $u_{n,x}(x,t), u_{n,xx}(x,t)$ vanish at infinities, $n=\overline{1,k+m}$.  

We are interested in a regular solution of problem (\ref{masala}). Moreover, the functions and their derivatives involved in the equation vanish as $x \to \infty$.

It should be noted that, if $\alpha_1=\alpha_2=...=\alpha_{k+m},$ then the Cauchy problem (\ref{masala}) has been solved in \cite{Sibir}. Moreover, if $\alpha_1=\alpha_2=...=\alpha_{k+m}=1,$ then the Airy-type equation coincides with the KdV equation with the standard derivative, and the Cauchy problem for this equation has been solved in \cite{Maqsad3}.

The Korteweg–de Vries (KdV) equation and its linearized Airy-type form play a fundamental role in the theory of dispersive wave propagation. These equations arise naturally in the study of shallow water waves, plasma physics, and transport processes in dispersive media. Initial-boundary value problems for third-order equations were studied in \cite{Himonas2019, pelloni2004}. These works highlight the mathematical and applied importance of dispersive evolution equations.

Parallel to this development, fractional differential equations (FDEs) have attracted considerable attention due to their ability to model memory effects and nonlocal phenomena more accurately than classical integer-order models. Their wide range of applications in physics, engineering, and complex systems is discussed in \cite{Mainardi1,  Hilfer1, Kilbas1, FM01, Metzler}. The fractional nature of such equations introduces hereditary properties into the dynamics, which are essential for describing anomalous diffusion and transport in heterogeneous media. 

Metric graphs provide a natural mathematical framework for modeling branched domains and complex networks in which transport, wave propagation, and diffusion processes occur. The theory of differential equations on geometric (metric) graphs has been systematically developed within the well-known Voronezh mathematical school founded by Yu. V. Pokornyi. This research direction focuses on boundary value problems, qualitative properties of solutions, spectral analysis, and applications to network-type physical models. Fundamental contributions of this school include the development of comparison principles, oscillation theory, and Green’s function techniques for differential operators on graphs, as presented in the monograph by Pokornyi, Penkin, Borovskikh, Pryadiev, Lazarev, and Shabrov \cite{pokornyi}. Early results on comparison theorems and boundary value problems for equations on graphs were obtained in works such as \cite{pokornyi-penkin}, which laid the groundwork for subsequent studies. Later investigations further extended this framework to spectral and inverse problems on spatial networks and metric graphs \cite{yurko}, confirming the continuing influence of the Voronezh school in modern graph-based differential equation theory. Such graphs arise in quantum graph theory, nanostructures, biological systems, and engineered network models \cite{Phys3, uzmat2022}. Fractional-order differential equations defined on graphs have been widely used to model processes with memory and hereditary effects. For example, in \cite{Ali}, a boundary value problem for a fractional differential equation with Caputo derivatives was studied on the cyclohexane graph representing the molecular structure of $C_6H_{12}$. In \cite{Mehandiratta}, the authors studied a space–time fractional parabolic equation on a metric star graph involving a Caputo time-fractional derivative and a spatial fractional operator of Sturm–Liouville type. Using the Galerkin method and functional analytic techniques, they proved the existence and uniqueness of a weak solution in appropriate fractional Sobolev spaces. In \cite{zhang}, the glycerol molecule was modeled as a graph, and existence of solutions was established under integral and derivative-dependent vertex conditions. These examples demonstrate that fractional dynamics on graphs provide an effective tool for describing complex transport mechanisms in branched media.

The combination of dispersive Airy-type equations and fractional dynamics on graphs leads to new analytical challenges. On metric graphs, a well-posed problem for an Airy-type equation was first formulated in \cite{Maqsad3}. Further studies of third-order dispersive equations on graphs include \cite{Calvante1, Noja1, Seifert1}. In \cite{Sibir}, the properties of potentials constructed from the fundamental solutions of the time-fractional Airy equation were investigated, and integral representations were derived for problems on infinite and finite intervals. In realistic settings, different channels of propagation may possess different internal structures or viscosities \cite{lobachevskiy2022}. Consequently, the temporal memory characteristics of the medium may vary from branch to branch, naturally leading to models with edge-dependent fractional orders. However, in most existing works the fractional order is assumed to be identical on all edges of the graph. This assumption simplifies the analysis but does not reflect heterogeneous media where memory effects vary across branches. Moreover, classical analytical tools, such as the Gronwall–Bellman inequality, are no longer applicable in the heterogeneous setting.

The present work is devoted to the Cauchy problem for a time-fractional Airy-type equation on an infinite star graph with edge-dependent Caputo derivatives. The heterogeneity of the parameters $\alpha_n$ models different temporal memory effects along different branches of the network. To formulate a mathematically consistent problem, we impose transmission conditions at the junction that ensure continuity of the solution and balance of fractional fluxes across the edges. To construct the solution for the problem we use the method of potentials. The vertex transmission conditions, which generalize Kirchhoff-type balance laws to the fractional framework, generate a coupled system that reduces to a generalized fractional Abel's integral equation. This reduction requires a refined analytical treatment. The solvability of the resulting generalized Abel's equation is established via results of A.~V.~Pskhu \cite{AVP03}, which allows us to prove existence and uniqueness of the solution and to obtain its explicit integral representation in terms of fractional Airy potentials. We note that an alternative approach to solving generalized Abel's integral equations can be found in \cite{li2021}, where inverse operator methods are employed. In contrast, our approach is based on the method of potentials. The developed framework provides a mathematical foundation for modeling dispersive transport processes with heterogeneous memory effects in branched network structures. It is worth noting that the problem under consideration can be formulated as a system of fractional partial differential equations (PDEs) with specific boundary conditions. Such systems have a wide range of applications and have been the subject of intensive investigation in recent years (see \cite{AshurovUmarov2022, Ashurov2025, Umarov2024} and references therein). Even in the context of systems of equations, the results presented here are novel.

\section{PRELIMINARIES}\label{asosiy-tushuncha}

In this section we provide the fundamental solutions to equation
\begin{equation}\label{teng-ma}\partial_{0t}^{\rho}u(x,t)-u_{xxx}(x,t)=f(x,t),\end{equation}
present some properties of potentials and formulate a theorem on the generalized Abel's integral equation from the work of A. Pskhu \cite{AVP03}.

The fractional integration operator is defined by formula \cite{Kilbas1}
\begin{equation} \label{integral}
D^{-\rho}_{0t}g(t) \equiv I^\rho_{0t}g(t) := \frac{1}{\Gamma(\rho)} \int\limits_0^{t} \frac{g(\xi)}{(t-\xi)^{1-\rho}} \, d\xi.
\end{equation}



For equation (\ref{teng-ma}), a fundamental solution was found in the form \cite{Pskhu2019}
\begin{equation} \label{fundamental1}
G_\rho^{2\rho/3}(x,t)=\frac{1}{3t^{1-2\rho/3}}\left\{
    \begin{array}{ll}
    \phi(-\rho/3,2\rho/3;\frac{x}{t^{\rho/3}}), &  x<0,\\
    -2\textrm{Re}[e^{{2\pi i}/3}\phi(-\rho/3,2\rho/3;e^{{2\pi i}/3}\frac{x}{t^{\rho/3}})],&  x>0.
    \end{array}\right.
\end{equation}

In \cite{Sibir}, the second elementary solution was found in the form
\begin{equation} \label{fundamental2}
V_\rho^{2\rho/3}(x,t) = \frac{1}{3t^{1-2\rho/3}} \textrm{Im}[e^{{2\pi i}/3} \phi(-\rho/3, 2\rho/3; e^{{2\pi i}/3} \frac{x}{t^{\rho/3}})], \quad x>0.
\end{equation}

Here, $\phi(\lambda,\mu;z)$ is the Wright function, defined by the relation \cite{FM01}
$$
\phi(\lambda,\mu;z) := \sum\limits_{n=0}^\infty \frac{z^n}{n! \, \Gamma(\lambda n + \mu)}, \quad \lambda > -1, \; \mu \in \textrm{C}.
$$

In the limit case $\rho\to 1$ the solutions \eqref{fundamental1} and  \eqref{fundamental2} get the form $G_1^{2/3}(x,t) = \frac{1}{(3t)^{1/3}}Ai\left(-\frac{x}{(3t)^{1/3}}\right)$,   
$V_1^{2/3}(x,t) = \frac{1}{(3t)^{1/3}}Bi\left(-\frac{x}{(3t)^{1/3}}\right)$, where $Ai$ and $Bi$ are Airy functions. We notice that these functions are fundamental and elementary solution to the Airy equation with $\rho=1$ (see \cite{Maqsad3}). 

The expression describing the exponentially decreasing behavior of the fundamental solution can be rewritten as follows (see \cite{Pskhu2019})

\begin{equation}\nonumber
\lim_{z\to \pm\infty}
t^{1-2\rho/3}G_{\rho}^{2\rho/3}(x,t)
\exp\!\left(
\nu_{\pm}x^{\frac{3}{3-\rho}}t^{-\frac{\rho}{3-\rho}}
\right)
=0,
\qquad
z=xt^{-\frac{\rho}{3}},
\end{equation}
where $ \nu_{-}<(3-\rho)\,3^{-\frac{3}{3-\rho}}\sigma^{\frac{\rho}{3-\rho}},
\nu_{+}=\nu_{-}\cos\frac{\pi}{3-\rho}.$ For the function $V_\rho^{2\rho/3}(x,t)$ the similar estimate is hold. We notice that if  $\rho=1$, the fundamental solution loses the exponentially decreasing property at $x\to +\infty.$

Based on solutions (\ref{fundamental1}) and (\ref{fundamental2}), in the work \cite{Sibir, rakhimov} were constructed the potentials and studied their properties. Let us mention them.

\begin{lemma} \cite{rakhimov}
\emph{Let functions} $\tau_k(t),\ (k=1,2)$ {be continuous in $[0,T]$. Then}

1. {Functions} $$w_1(x,t)=\int_0^t{G_{\rho}^{2\rho/3}(x-a,t-\eta )\tau_1}(\eta )d\eta$$ {and} $$w_2(x,t)=\int_0^t V_\rho^{2\rho/3}(x-a,t-\eta )\tau_2(\eta )d\eta$$  {are solutions of the equation} $$\partial_{0t}^\rho u_j(x,t)-\frac{\partial^3 u_j(x,t)}{\partial x^3}=0;$$

2. {For the functions} $w_1(x,t)$ {and} $w_2(x,t)$ {hold relations} $$\lim_{t\to 0} w_k(x,t)=0, k=1,2.$$

\end{lemma}

\begin{lemma} \cite{rakhimov} {We put}
$$w_3(x,t)=\int_0^t\frac{\partial^2}{\partial x^2}G_{\rho}^{2\rho/3}(x-a,t-\eta)\tau_3(\eta) d\eta $$ {and}
$$w_4(x,t)=\int_0^t\frac{\partial^2}{\partial x^2}V_\rho^{2\rho/3}(x-a,t-\eta )\tau_4(\eta) d\eta.$$
{Let} $\tau_3(t), \tau_4(t)$ be continuous  on $[0,T].$ {Then}
$$\lim_{x\to a-0} w_3(x,t)=\frac{1}{3}\tau_3(t),$$
$$\lim_{x\to a+0} w_3(x,t)=-\frac{2}{3}\tau_3(t),$$
$$\lim_{x\to a+0} w_4(x,t)=0.$$

\end{lemma}

\begin{lemma} \cite{rakhimov} 
{Let} $\tau_5(x) \in C[a,b].$ {Then the function $$w_5(x,t)=\int_a^b G_\rho^{2\rho/3}(x-\xi,t )\tau_5(\xi)d\xi$$ is a solution for equation $$\partial_{0t}^{\rho}u(x,t)-u_{xxx}(x,t)=f(x,t)$$ and}
$$
w_5(x,0)=\tau_5(x).
$$

\end{lemma}

\begin{lemma}  \cite{rakhimov} 
{Let $f(x,t)\in C^{0,1}([a,b]\times[0,T]).$ The equation} $\partial_{0t}^{\rho}u(x,t)-u_{xxx}(x,t)=f(x,t)$ {with initial condition}
$$u(x,0)=0$$ {has a solution in the form}
$$
w_6(x,t)=\int_{0}^{t}d\eta\int_{a}^{b}G_{\rho}^{2\rho/3}(x-\xi,t-\eta)f(\xi,\eta)d\xi.
$$
\end{lemma}

The proofs of these lemmas are given in the work  \cite{rakhimov}.

For our further investigation we need to give solution to the generalized Abel's equation
\begin{equation}\label{abel}
    \sum\limits_{i=0}^{n}A_iI_{0t}^{\rho_i}x(t)=g(t),
\end{equation}
where $A_0\neq 0,\,\,\, \rho_i\geq 0,$ for $i=\overline{1,n}.$
Following to \cite{AVP03}, we introduce the function
\begin{equation}\label{S-funksiya}
    S_n^\nu(t;z_1,...,z_n;\gamma_1,...,\gamma_n)=(h_1*h_2*...*h_n)(t),
\end{equation}
where $(h_i*h_j)(t)=\int\limits_0^th_i(t-\tau)h_j(\tau)d\tau$ and
$h_i(t)=t^{\nu_i-1}\phi(\gamma_i,\nu_i;z_it^{\gamma_i}).$
The parameters and arguments of the function (\ref{S-funksiya}) satisfy the conditions
$$t>0,\,\,\,z_i\in\mathbf{R}, \,\,\, \gamma_i>0, \,\,\, \nu_i>0, \,\,\, \nu=\sum\limits_{i=1}^n\nu_i.$$
Using the function (\ref{S-funksiya}), we define another function
$$G_n^\nu(t;\gamma_1,...,\gamma_n;\gamma_1,...,\gamma_n)=\int\limits_0^{+\infty}e^{-\tau}S_n^\nu(\tau;\gamma_1\tau,...,\gamma_n\tau;\gamma_1,...,\gamma_n)=(h_1*h_2*...*h_n)(t)dt.$$
Note that this function should not be confused with the function (\ref{fundamental1}).

\begin{theorem}\label{pskhu-theo} \cite{AVP03} Let $A_0\neq 0,$ $0\leq\rho_0<\rho_k,$ $g(t)\in L(0,T).$ The equation (\ref{abel}) has a unique integrable solution $x(t),$ which for each $\nu>\rho_0$ can be represented in the form
\begin{equation} \label{abel-yechim}
    x(t)=\partial_{0t}^\nu\left( g*\omega_\nu\right)(t),
\end{equation}
where $\omega_\nu=\omega_\nu(t)$ is given by
$$\omega_\nu(t)=G_n^{\nu}\left(t;-\frac{A_1}{A_0};...,-\frac{A_m}{A_0};\rho_1-\rho_0,...,\rho_n-\rho_0 \right).$$
\end{theorem}

Chenkuan Li and Hari M. Srivastava \cite{li2021} conducted a rigorous analysis of the generalized Abel's equation, presenting its solution in the form of a power series. This construction relies on the resolvent operator, Mittag-Leffler functions, and the Laplace transform. Such a solution is a significant and practical tool for addressing applied problems. Additionally, their work establishes the unique solvability of a system of nonlinear generalized Abel's integral equations, further enhancing its theoretical and practical relevance.

\section{UNIQUENESS OF THE SOLUTION}

In this section we prove the uniqueness of the solution to the problem (\ref{masala}).

\begin{theorem} \label{Uniqueness} Let $B^TB-I$ be a negatively definite matrix. Then problem (\ref{masala}) has at most one solution.
\end{theorem}
{\sf Proof.}
Let it have two different solutions $u^{(1)}=(u^{(1)}_1, u^{(1)}_2,...,u^{(1)}_{k+m})^T$ and $u^{(2)}=(u^{(2)}_1, u^{(2)}_2,...,u^{(2)}_{k+m})^T.$ Then the function $v=u^{(1)}-u^{(2)}$ satisfies the equation

$$\partial_{0t}^{\alpha_j}v_j(x,t)=v_{j,xxx}(x,t),$$
homogeneous initial condition $v_j(x,0)=0,$ and 3-5 conditions of the system (\ref{masala}).

We multiply both sides of the equation $\partial_{0t}^{\alpha_n} v_n(x,t) = v_{n,xxx}(x,t)$ by $v_n(x,t)$, use the inequality \cite{Ali1}
$$
v \partial_{0t}^{\alpha} v \geq \frac{1}{2} \partial_{0t}^{\alpha} v^2,
$$
and integrate the resulting expression over $B_j$. Taking into account the inequality
$$
\int\limits_{a}^b \partial_{0t}^{\alpha_n} v_n^2 \, dx \leq 2 \int\limits_{a}^b v_n \partial_{0t}^{\alpha_n} v_n \, dx = 2 \int\limits_{a}^b v_n v_{n,xxx} \, dx = 2 \left( \left. v_n v_{n,xx} \right|_{a}^b - \frac{1}{2} \left. v_{n,x}^2 \right|_{a}^b \right),
$$
summing up such inequalities over all edges, we obtain
$$
\sum\limits_{n=1}^{k+m} \partial_{0t}^{\alpha_n} \int\limits_{b_n} v_n^2 \, dx \leq 2 \sum\limits_{i=1}^{k} \left. \left( v_i v_{i,xx} \right) \right|_{-\infty}^{0} + 2 \sum\limits_{j=k+1}^{k+m} \left. \left( v_j v_{j,xx} \right) \right|_{0}^{+\infty} + \sum\limits_{i=1}^{k} v_{i,x}^2(0,t) - \sum\limits_{j=k+1}^{k+m} v_{j,x}^2(0,t).
$$
From this inequality, we have
$$
\sum_{n=1}^{k+m} \partial_{0t}^{\alpha_n} \int_{b_n} v_n^2 \, dx
\le
2v_1(0,t)\sum\limits_{i=1}^{k}\frac{1}{a_i^2}v_{i,xx}(0,t)
- 2v_1(0,t) \sum_{j=k+1}^{k+m} \frac{1}{a_j^2} v_{j,xx}(0,t)+(v^-(0,t))^Tv^-(0,t)- $$
$$-(v^+(0,t))^Tv^+(0,t)=2v_1(0,t) \left( \sum\limits_{i=1}^{k}\frac{1}{b_i^2}v_{i,xx}(0,t)
-  \sum_{j=k+1}^{k+m} \frac{1}{b_j^2} v_{j,xx}(0,t) \right)-$$
$$+(v^-(0,t))^Tv^-(0,t)-(v^-(0,t))^TB^TBv^-(0,t).$$
Taking into account the conditions of the theorem and vertex conditions \eqref{masala}, we obtain
\[
\sum\limits_{n=1}^{k+m} \partial_{0t}^{\alpha_n} \int\limits_{b_n} v_n^2 \, dx \leq \left(v^-(0,t)\right)^T \left( B^T B - I \right) v^-(0,t) \leq 0.
\]
Integrating this inequality and using the homogeneous initial condition, we have
\[
 \sum\limits_{n=1}^{k+m}  I_{0t}^{1 - \alpha_n} \int\limits_{b_n} v_n^2 \, dx  \leq 0.
\]

In the last inequality, the sum of several nonnegative terms is less than or equal to zero. Therefore, each term must be zero, which implies $v \equiv 0$. Hence, $u^{(1)} \equiv u^{(2)}$. The theorem is proved.

The condition $B^{T}B - I < 0$ has a clear physical and geometric interpretation. From a physical point of view, it means that the vertex does not generate energy. Indeed, the matrix $B$ describes how spatial fluxes (or derivatives) incoming to the vertex are redistributed among the outgoing edges. The inequality $||Bv||<||v||$ implies that the total outgoing flux is strictly smaller than the incoming one, which corresponds to a dissipative or passive junction. As a consequence, no amplification of waves occurs at the vertex, ensuring physical stability of the model.
From a geometric perspective, the matrix $B$ acts as a contraction operator in the corresponding Euclidean space. This contraction property guarantees that the energy functional associated with the homogeneous problem is non-increasing in time. Such a property is crucial for proving uniqueness, since it prevents the existence of nontrivial solutions with zero initial data.

\textbf{An Illustrative Example 1.}
Here we give simple example which demonstrate that if $B^TB-I$ is not a negatively definite then the solution of the problem
can be more than one. 

Consider problem \eqref{masala} on a star graph with one incoming and one outgoing edge, i.e., $k=m=1$. 
Let
\[
\alpha_1=\alpha_2=\alpha\in(0,1),\qquad a_2=1,\qquad B=(b).
\]
We choose $b=-2$.
Then
$B^TB-I=b^2-1=4-1=3>0,$
so the condition of Theorem 3.1 is violated.

Assume 
$
f_1(x,t)=f_2(x,t)=0, u_{0,1}(x)=u_{0,2}(x)=0.
$
The problem becomes
\[
\begin{cases}
\partial_{0t}^{\alpha}u_1(x,t)-u_{1,xxx}(x,t)=0, & x<0,\ t>0,\\[1mm]
\partial_{0t}^{\alpha}u_2(x,t)-u_{2,xxx}(x,t)=0, & x>0,\ t>0,\\[1mm]
u_1(0,t)=u_2(0,t),\\[1mm]
u_{2,x}(0,t)=-2u_{1,x}(0,t),\\[1mm]
u_{1,xx}(0,t)=u_{2,xx}(0,t).
\end{cases}
\]

In addition to the trivial solution, this problem has a non-trivial solution, which is given by the formula
 $$u_1(x,t)=t \phi\left(-\frac{\alpha}{3},2,\frac{x}{t^\frac{\alpha}{3}}\right)$$ and
$$u_2(x,t)=2t \textrm{Re} \left[\phi\left(-\frac{\alpha}{3},2,\frac{\omega x}{t^\frac{\alpha}{3}}\right)\right],
$$
where
 $\omega=-\frac{1}{2}+i\frac{\sqrt3}{2}$. 
Therefore, the solution is not unique.

\textbf{An Illustrative Example 2.}
In this example we want to demonstrate the condition of the Lemma \ref{Uniqueness} in the cases $k=1, m=2$ and $k=2, m=1$. 

\noindent
\textbf{Case 1:} \(k=1,\ m=2\). In this case, the matrix \(B\) has size \(2\times 1\), so we write
$B=
\begin{pmatrix}
b_1\\
b_2
\end{pmatrix}.$
Then
$$B^{T}B-I=
\begin{pmatrix}
b_1 & b_2
\end{pmatrix}
\begin{pmatrix}
b_1\\
b_2
\end{pmatrix}-I
=
(b_1^2+b_2^2-1).
$$

Therefore, the condition of Lemma 3.1 is satisfied whenever
$b_1^2+b_2^2<1.$

\medskip

\noindent
\textbf{Case 2:} \(k=2,\ m=1\). In this case, the matrix \(B\) has size \(1\times 2\), so we write
\[
B=
\begin{pmatrix}
b_1 & b_2
\end{pmatrix}.
\]
Then
\[
B^{T}B-I=
\begin{pmatrix}
b_1\\
b_2
\end{pmatrix}
\begin{pmatrix}
b_1 & b_2
\end{pmatrix}
-I=
\begin{pmatrix}
b_1^2-1 & b_1b_2\\
b_1b_2 & b_2^2-1
\end{pmatrix}.\]
So, using Sylvester's criterion we conclude that $B^TB-I$ is negatively defined iff $b_1^2+b_2^2<1$.

\section{EXISTENCE OF THE SOLUTION}\label{mavjudlik}

In the works\cite{Sibir} the potential theory for time-fractional Airy-type equations with fractional derivatives was developed. Based on this theory, we will seek the solution of the problem in the following form

\begin{equation} \label{yechim}
u_n(x,t) = \int\limits_{0}^{t} G_{\alpha_n}^{2\alpha_n/3}(x-0, t-\tau) \lambda_n(\tau)  d\tau + \int\limits_{0}^{t} V_{\alpha_n}^{2\alpha_n/3}(x-0, t-\tau) \mu_n(\tau)  d\tau + F_n(x,t), \quad n = \overline{1,k+m},
\end{equation}
where
\[
F_n(x,t) = \int\limits_{b_n} G_{\alpha_n}^{2\alpha_n/3}(x-\xi, t-0) u_{0,n}(\xi)  d\xi + \int\limits_{0}^{t} \int\limits_{b_n} G_{\alpha_n}^{2\alpha_n/3}(x-\xi, t-\tau) f_n(\xi,\tau)  d\xi  d\tau.
\]

According to the lemmas 1--4, the functions \eqref{yechim} satisfy the equation and initial condition in \eqref{masala}. 
Here, the functions $\lambda_n(t)$ and $\mu_n(t)$ are unknowns and must be determined. To do this, we will use the remaining vertex conditions from the \eqref{masala}. As the number of unknowns are more than the number of the vertex conditions, we put  $\mu_1(t) =\mu_2(t) =...=\mu_k(t) = 0$ to equalize them. 

Throughout the analysis, we assume without loss of generality that $\alpha_1=\max\{\alpha_1,...,\alpha_{k+m}\}.$
This assumption is purely technical and does not restrict the generality of the results. Indeed, the labeling of the edges of the star graph is arbitrary, and if for some index $j$ one has $\alpha_j > \alpha_1,$ the same analysis can be carried out after a simple relabeling of the edges, taking this index as the reference one. The above choice is made only to simplify the representation of the arising fractional integral operators.

First, we consider the first vertex condition $Au(0,t)=0$ in \eqref{masala}.

\begin{align*}
a_j & \left( \int\limits_{0}^{t} G_{\alpha_j}^{2\alpha_j/3}(0, t-\tau) \lambda_j(\tau)  d\tau + \int\limits_{0}^{t} V_{\alpha_j}^{2\alpha_j/3}(0, t-\tau) \mu_j(\tau)  d\tau \right) \\
&= \int\limits_{0}^{t} G_{\alpha_1}^{2\alpha_1/3}(0, t-\tau) \lambda_1(\tau)  d\tau d\tau + F_1(0,t) - a_j F_j(0,t),
\end{align*}
where $j = \overline{k+1, k+m}$.

Substituting the values of $G(0,t)$ and $V(0,t)$, we obtain
\[
\frac{a_j}{3\Gamma\left( \frac{2\alpha_j}{3} \right)} \int\limits_{0}^{t} \frac{\lambda_j(\tau) + \frac{\sqrt{3}}{2} \mu_j(\tau)}{(t - \tau)^{1 - 2\alpha_j/3}}  d\tau = \frac{1}{3\Gamma\left( \frac{2\alpha_1}{3} \right)} \int\limits_{0}^{t} \frac{\lambda_1(\tau)}{(t - \tau)^{1 - 2\alpha_1/3}}  d\tau + F_1(0,t) - a_j F_j(0,t),
\]
for $j = \overline{2, k+m}$.

Using the definition of the fractional integral, we get
\[
a_j \left( \lambda_j(t) + \frac{\sqrt{3}}{2} \mu_j(t) \right) = I_{0t}^{\frac{2(\alpha_1 - \alpha_j)}{3}} \lambda_1(t) + 3 \, D_{0t}^{\frac{2\alpha_j}{3}} \left( F_1(0,t) - a_j F_j(0,t) \right), \quad j = \overline{2, k+m}.
\]

From this, it follows that
\begin{equation} \nonumber
    {{\lambda }_{j}}\left( t \right)+\frac{\sqrt{3}}{2}{{\mu }_{j}}\left( t \right)=\frac{1}{{{a}_{j}}}I_{0t}^{\frac{2\left( {{\alpha }_{1}}-{{\alpha }_{j}} \right)}{3}}{{\lambda }_{1}}\left( t  \right)+\frac{3}{{{a}_{j}}}D _{0t}^{\frac{2{{\alpha }_{j}}}{3}}\left( {{F}_{1}}\left( 0,t \right)-{{a}_{j}}{{F}_{j}}\left( 0,t \right) \right),\,\,\,\,j=\overline{2,k+m}.
\end{equation}

Taking into account that $\mu_1 = \mu_2 = \dots = \mu_k = 0$, from the last relation we find

\begin{equation} \label{c1r1}
{{\lambda }_{i}}\left( t \right)=\frac{1}{{{a}_{i}}}I_{0t}^{\frac{2\left( {{\alpha }_{1}}-{{\alpha }_{i}} \right)}{3}}{{\lambda }_{1}}\left( t \right)+\frac{3}{{{a}_{i}}}D_{0t}^{\frac{2{{\alpha }_{i}}}{3}}\left( {{F}_{1}}\left( 0,t \right)-a_i{{F}_{i}}\left( 0,t \right) \right),\,\,\,\,i=\overline{2,k},
\end{equation}
\begin{equation} \label{c1r2}
{{\lambda }_{j}}\left( t \right)+\frac{\sqrt{3}}{2}{{\mu }_{j}}\left( t \right)=\frac{1}{{{a}_{j}}}I_{0t}^{\frac{2\left( {{\alpha }_{1}}-{{\alpha }_{j}} \right)}{3}}{{\lambda }_{1}}\left( t \right)+\frac{3}{{{a}_{j}}}D_{0t}^{\frac{2{{\alpha }_{j}}}{3}}\left( {{F}_{1}}\left( 0,t \right)-a_j{{F}_{j}}\left( 0,t \right) \right),\,\,\,\,j=\overline{k+1,k+m}.
\end{equation}

Next, we compute the derivative of the function
\[
u_{j,x}(x,t) = \int\limits_{0}^{t} G_{\alpha_j}^{\alpha_j/3}(x-0, t-\tau) \lambda_j(\tau)  d\tau + \int\limits_{0}^{t} V_{\alpha_j}^{2\alpha_j/3}(x-0, t-\tau) \mu_j(\tau)  d\tau + F_{j,x}(x,t),
\]
and evaluate it at $x = 0$
\[
u_{j,x}(0,t) = \frac{1}{3\Gamma\left( \frac{\alpha_j}{3} \right)} \int\limits_{0}^{t} \frac{\lambda_j(\tau)}{(t - \tau)^{1 - \alpha_j/3}}  d\tau - \frac{\sqrt{3}}{6\Gamma\left( \frac{\alpha_j}{3} \right)} \int\limits_{0}^{t} \frac{\mu_j(\tau)}{(t - \tau)^{1 - \alpha_j/3}}  d\tau + F_{j,x}(0,t).
\]

We rewrite the second vertex condition $u_x^+(0,t)=Bu_x^-(0,t)$ in \eqref{masala} in the form
\[
u_{j,x}(0,t) = \sum\limits_{i=1}^{k} b_{ji} u_{i,x}(0,t), \quad j = \overline{k+1, k+m},
\]
and substitute the expressions derived above to obtain

\begin{equation}\label{c2r}
I_{0t}^{\alpha_j/3} \lambda_j(t) - \frac{\sqrt{3}}{2} I_{0t}^{\alpha_j/3} \mu_j(t) = \sum\limits_{i=1}^{k} b_{ji} \left( I_{0t}^{\alpha_i/3} \lambda_i(t) + 3 F_{i,x}(0,t) \right) - 3 F_{j,x}(0,t).
\end{equation}

Substituting \eqref{c2r} into \eqref{c1r1}, we obtain
\begin{equation}\label{sist-1}
\lambda_j(t) - \frac{\sqrt{3}}{2} \mu_j(t) =
b_{j1} I_{0t}^{\frac{\alpha_1 - \alpha_j}{3}} \lambda_1(t) +
\sum_{i=2}^{k} \frac{b_{ji}}{a_i} I_{0t}^{\frac{2\alpha_1 - \alpha_i - \alpha_j}{3}} \lambda_1(t) +
g_j(t), \quad j = \overline{k+1, k+m}.
\end{equation}
where
\begin{align*}
g_j(t) &= 3 \sum_{i=2}^{k} \frac{b_{ji}}{a_i} D_{0t}^{\alpha_j/3} I_{0t}^{\alpha_i/3} D_{0t}^{2\alpha_i/3} \left( F_1(0,t) - a_i F_i(0,t) \right) \\
&\quad + 3 \sum_{i=1}^{k} b_{ji} D_{0t}^{\alpha_j/3} F_{i,x}(0,t) - 3 D_{0t}^{\alpha_j/3} F_{j,x}(0,t), \qquad j = \overline{k+1, k+m}.
\end{align*}

From equations \eqref{c1r2} and \eqref{sist-1} we find
\begin{equation}\label{lambda_j}
\lambda_j(t) = \frac{1}{2a_j} I_{0t}^{\frac{2(\alpha_1 - \alpha_j)}{3}} \lambda_1(t) + \frac{b_{j1}}{2} I_{0t}^{\frac{\alpha_1 - \alpha_j}{3}} \lambda_1(t) + \sum\limits_{i=2}^{k} \frac{b_{ji}}{2a_i} I_{0t}^{\frac{2\alpha_1 - \alpha_i - \alpha_j}{3}} \lambda_1(t) + \frac{h_j(t) + g_j(t)}{2}
\end{equation}
and
\begin{equation}\label{mu_j}
\mu_j(t) = \frac{1}{\sqrt{3}a_j} I_{0t}^{\frac{2(\alpha_1 - \alpha_j)}{3}} \lambda_1(t) - \frac{b_{j1}}{\sqrt{3}} I_{0t}^{\frac{\alpha_1 - \alpha_j}{3}} \lambda_1(t) - \sum\limits_{i=2}^{k} \frac{b_{ji}}{\sqrt{3}a_i} I_{0t}^{\frac{2\alpha_1 - \alpha_i - \alpha_j}{3}} \lambda_1(t) + \frac{h_j(t) - g_j(t)}{\sqrt{3}},
\end{equation}
where $j = \overline{k+1, k+m}$, and the functions $h_j(t)$ are defined by
\begin{align*}
h_j(t) &= \frac{3}{a_j} D_{0t}^{\frac{2\alpha_j}{3}} \left( F_1(0,t) - a_j F_j(0,t) \right).
\end{align*}

Next, using Lemma 3, we find the values of the second derivatives at $x = 0$
\[
\lim_{x \to -0} u_{1,xx}(x,t) = \frac{1}{3} \lambda_1(t) + F_{1,xx}(0,t),
\]
and
\[
\lim_{x \to +0} u_{j,xx}(x,t) = -\frac{2}{3} \lambda_j(t) + F_{j,xx}(0,t), \quad j = \overline{2, k+m}.
\]

Substituting these into the last vertex condition in \eqref{masala}, we obtain
\begin{equation}\label{c3r}
\lambda_1(t) + 3F_{1,xx}(0,t) + \sum\limits_{i=2}^{k} \frac{1}{a_i} \left( \lambda_i(t) + 3F_{i,xx}(0,t) \right) =
\sum\limits_{j=k+1}^{k+m} \frac{1}{a_j} \left( -2\lambda_j(t) + 3F_{j,xx}(0,t) \right).
\end{equation}
    
    Now substitute relations (\ref{c1r2}) and (\ref{lambda_j}) into (\ref{c3r}). After simplification, we get the following generalized Abel's integral equation in the form
    \begin{align}\label{lambda1}
    \lambda_1(t) &+ \sum\limits_{i=2}^{k} \frac{1}{a_i^2} I_{0t}^{\frac{2(\alpha_1 - \alpha_i)}{3}} \lambda_1(t)
    + \sum\limits_{j=k+1}^{k+m} \frac{1}{a_j^2} I_{0t}^{\frac{2(\alpha_1 - \alpha_j)}{3}} \lambda_1(t)
    + \sum\limits_{j=k+1}^{k+m} \frac{b_{j1}}{a_j} I_{0t}^{\frac{\alpha_1 - \alpha_j}{3}} \lambda_1(t)
    \nonumber \\
    &+ \sum\limits_{j=k+1}^{k+m} \sum\limits_{i=2}^{k} \frac{b_{ji}}{a_i a_j} I_{0t}^{\frac{2\alpha_1 - \alpha_i - \alpha_j}{3}} \lambda_1(t)
    =
    \sum\limits_{j=k+1}^{k+m} \left( \frac{1}{a_j} F_{j,xx}(0,t) - \frac{h_j(t) + g_j(t)}{a_j} \right)
    \\
    &- F_{1,xx}(0,t) - \sum\limits_{i=2}^{k} \left( \frac{3}{a_i^2} D_{0t}^{\frac{2\alpha_i}{3}} \left( F_1(0,t) - a_i F_i(0,t) \right) + \frac{1}{a_i} F_{i,xx}(0,t) \right).
    \nonumber
    \end{align}

Its solution is presented in Section~2. Let us rewrite this equation in the form
\begin{equation}\label{abel1}
\sum_{s=0}^{km+k+m} A_s I_{0t}^{\rho_s}\lambda_1(t)=g(t),
\end{equation}
where $A_0=1, \qquad \rho_0=0,\qquad$ $A_{i-1}=\frac{1}{a_i^2},$
$\qquad
\rho_{i-1}=\frac{2(\alpha_1-\alpha_i)}{3},
\qquad (i=\overline{2,k}),$ $$A_{j-1}=\frac{1}{a_j^2},
\qquad \rho_{j-1}=\frac{2(\alpha_1-\alpha_j)}{3},
\qquad (j=\overline{k+1,k+m}),$$ 

$A_{m+j-1}=\frac{b_{j1}}{a_j},
\qquad \rho_{m+j-1}=\frac{\alpha_1-\alpha_j}{3},
\qquad (j=\overline{k+1,k+m}),
$
$$A_{(k-1)+2m+(j-k-1)(k-1)+(i-1)} 
=\frac{b_{ji}}{a_i a_j},$$
 $\rho_{(k-1)+2m+(j-k-1)(k-1)+(i-1)}
=\frac{2\alpha_1-\alpha_i-\alpha_j}{3},
$
$ 
(j=\overline{k+1,k+m}, 
\quad 
i=\overline{2,k})
$ and
\begin{align}
g(t)=&
\sum_{j=k+1}^{k+m}
\left(
\frac{1}{a_j} F_{j,xx}(0,t)
-\frac{h_j(t)+g_j(t)}{a_j}
\right)-\nonumber
\\
&-\sum_{i=2}^{k}
\left(
\frac{3}{a_i^2} D_{0t}^{\frac{2\alpha_i}{3}}
\bigl( F_1(0,t)-a_iF_i(0,t) \bigr)
+\frac{1}{a_i} F_{i,xx}(0,t)
\right)- F_{1,xx}(0,t). \label{g(t)}
\end{align}

Relying on Theorem \ref{pskhu-theo}, we can find
\begin{equation} \label{oxir_yechim}
    \lambda_1(t)=\partial_{0t}^\nu\left( g*\omega_\nu\right)(t),
\end{equation}
where $\omega_\nu=\omega_\nu(t)$ is given by
$$\omega_\nu(t)=G_n^{\nu}\left(t;-A_1;-A_2;...,-A_m;\rho_1,...,\rho_n \right).$$

It is known that equation (\ref{abel1}) is of Volterra type, and to determine the class of $\lambda_1(t)$, it suffices to impose a condition on $g(t)$; namely, $g(t)\in C^1(0,T)$. Taking into account the expressions $F_{j}(0,t)$, $F_{j,x}(0,t)$, and $F_{j,xx}(0,t)$ appearing in $g(t)$, we can determine the function classes of $u_{0,n}(x)$ and $f_n(x,t)$. Summarizing, we formulate the following theorem.

\begin{theorem}  Let $B^T B - I$ be a negative definite matrix and the functions $f_n(x,t) \in C^{2,1}(\overline{b_n} \times [0,T])$ and $u_{0,n}(x) \in C^2(\overline{b_n})$ and their derivatives tend to zero at infinities. Then the problem~(\ref{masala}) has a unique solution on the form (\ref{yechim}). Here $\lambda_1(t)$ is is given by (\ref{oxir_yechim}), $\lambda_j(t)$ and $\mu_j(t)$ for $j = \overline{2, k+m}$ will be found by ~(\ref{lambda_j}) and~(\ref{mu_j}).
\end{theorem}

\textbf{An Illustrative Example 3.}
Consider problem \eqref{masala} on a star graph with one incoming and two outgoing edges, that is, \(k=1\), \(m=2\). Let
$\alpha_1=\frac{9}{10},\qquad \alpha_2=\frac{2}{3},\qquad \alpha_3=\frac{1}{2},
\qquad 
a_2=a_3=1,$
and
$B=\begin{pmatrix}\frac12, \frac13\end{pmatrix}^T$. For simplicity of presentation we put $f_n(x,t)=e^{-|x|},$ $u_{0.n}(x)=0$, $n=1,2,3$.

Then the main equation (i.e. generalized Abel's integral equation) \eqref{abel} has the form  
\[
\lambda_1(t)
+ I_{0t}^{\frac{7}{45}}\lambda_1(t)
+ I_{0t}^{\frac{4}{15}}\lambda_1(t)
+ \frac12 I_{0t}^{\frac{7}{90}}\lambda_1(t)
+ \frac13 I_{0t}^{\frac{2}{15}}\lambda_1(t)
+ \frac16 I_{0t}^{\frac{19}{90}}\lambda_1(t)
= g(t),
\]
where \(g(t)\) is determined by \eqref{g(t)}. By Theorem 2.1 it has a unique solution
\[
\lambda_1(t)=\frac{d}{dt}(g*\omega_1)(t),
\]
where
\[
\omega_1(t)=
G_5^1\!\left(
t;
-1,-1,-\frac12,-\frac13,-\frac16;
\frac{7}{45},\frac{4}{15},\frac{7}{90},\frac{2}{15},\frac{19}{90}
\right).
\]

The remaining densities \(\lambda_2,\lambda_3,\mu_2,\mu_3\) can be found from equation \eqref{c1r1}, \eqref{lambda_j} and \eqref{mu_j}. Substituting them into the potential representation, we obtain
\[
u_1(x,t)
=
\int_0^t G_{9/10}^{3/5}(x,t-\tau)\lambda_1(\tau)\,d\tau
+
\int_{b_1}G_{9/10}^{3/5}(x-\xi,t)e^{-|\xi|}\,d\xi,
\]
\[
\begin{aligned}
u_2(x,t)=\;&
\int_0^t G_{2/3}^{4/9}(x,t-\tau)\lambda_2(\tau)\,d\tau
+
\int_0^t V_{2/3}^{4/9}(x,t-\tau)\mu_2(\tau)\,d\tau \\
&+
\int_{b_2}G_{2/3}^{4/9}(x-\xi,t)e^{-|\xi|}\,d\xi,
\end{aligned}
\]
and
\[
\begin{aligned}
u_3(x,t)=\;&
\int_0^t G_{1/2}^{1/3}(x,t-\tau)\lambda_3(\tau)\,d\tau
+
\int_0^t V_{1/2}^{1/3}(x,t-\tau)\mu_3(\tau)\,d\tau \\
&+
\int_{b_3}G_{1/2}^{1/3}(x-\xi,t)e^{-|\xi|}\,d\xi.
\end{aligned}
\]
Thus, the solution is represented explicitly via the solution of the associated generalized Abel's integral equation.

\section{CONCLUSION}

In this paper, we studied the Cauchy problem for a time-fractional Airy-type equation on an infinite star graph with edge-dependent Caputo derivatives. The heterogeneity of the fractional orders models different temporal memory characteristics along the branches of the network, which is essential for describing dispersive transport processes in non-uniform media. By imposing transmission conditions that generalize Kirchhoff-type balance laws to the fractional setting, we formulated a mathematically consistent problem on the graph. Using the method of potentials constructed from the fundamental solutions of the fractional Airy operator, we reduced the problem to a system of Volterra-type integral equations in time.

The key analytical step consists in the reduction to a generalized fractional Abel's integral equation, whose solvability follows from results of A.~V.~Pskhu. This approach allowed us to establish existence and uniqueness of the solution and to obtain its explicit integral representation in terms of fractional Airy potentials. The developed framework provides a rigorous mathematical basis for the analysis of dispersive evolution equations with heterogeneous memory effects on branched network structures and opens the way for further investigations of more general graph configurations and nonlinear fractional models.

\end{document}